\documentclass{amsart}

\usepackage{amsmath}
\usepackage{amssymb}

\newtheorem{theorem}{Theorem}
\newtheorem{proposition}[theorem]{Proposition}

\theoremstyle{definition}

\theoremstyle{remark}

\begin{document}
\title{A note on the geometry of figurate numbers}
\author{Franti\v sek Marko}
\address{Pennsylvania State University \\ 76 University Drive \\ Hazleton, PA 18202 \\ USA}
\email{fxm13@psu.edu}
\begin{abstract}
We give a short proof of the formula 
$n^p=\sum_{\ell=0}^{p-1} (-1)^{\ell} c_{p,\ell} F^{p-\ell}_n$, where 
$F^{p-\ell}_n$ is the figurate number of dimension $p-\ell$, and $c_{p,\ell}$ is the number of $(p-\ell)$-dimensional facets  of 
$p$-dimensional simplices obtained by cutting the $p$-dimensional cube.
This formula was formulated as Conjecture 16 in \cite{ml}.
\end{abstract}
\maketitle
Recall from \cite{ml} that the $k$-dimensional cube $0\leq x_1, \ldots, x_k\leq n-1$ is cut into $n!$ $k$-dimensional simplices given by $x_{\sigma_1}\geq x_{\sigma_2}\geq \ldots \geq x_{\sigma_k}$ for a permutation $\sigma$ of the set $\{1, \ldots, k\}$. The $(k-\ell)$-dimensional facets of these simplices are given by $0\leq x_1, \ldots, x_k\leq n-1$ and the conditions $x_{\sigma_1}L_1x_{\sigma_2}L_2\ldots L_{k-1}x_{\sigma_k}$, where exactly $\ell$ symbols $L_i$ equal $``="$, and $k-1-\ell$ symbols $L_i$ equal $``\geq"$. Denote by $c_{k,\ell}$ the number of all such $(k-\ell)$-dimensional facets.

Also, denote bv $F^k_n=\binom{n+k-1}{n}$ the figure number of dimension $k$. In Section 5 of \cite{ml}, we gave a geometric counting argument using the inclusion-exclusion principle in support of Conjecture 16 of \cite{ml} stating the following.

\begin{proposition} For all positive integers $n$ and $p$, 
$n^p=\sum_{\ell=0}^{p-1} (-1)^{\ell} c_{p,\ell} F^{p-\ell}_n$.
\end{proposition}
\textit{Proof:}
First we observe that $c_{p,\ell}$ counts the number of surjective maps from the set $\{1, \ldots, p\}$ to the set 
$\{1, \ldots, p-\ell\}$. Indeed, the $p-1-\ell$ symbols $L_{a_1}, L_{a_2}, \ldots,$ $ L_{a_{p-1-\ell}}$ that are equal to $``\geq"$ in the expression $x_{\sigma_1}L_1x_{\sigma_2}L_2\ldots L_{p-1}x_{\sigma_p}$ separate indices $\{1, \ldots, p\}$ into $p-\ell$ nonempty groups $\{\sigma_1, \ldots, \sigma_{a_1}\}$, $\{\sigma_{a_1+1}, \ldots, \sigma_{a_2}\}$, 
$\{\sigma_{a_{p-\ell-2}+1}, \ldots, \sigma_{a_{p-1-\ell}}\}$, and 
$\{\sigma_{a_{p-1-\ell}+1}, \ldots \sigma_p\}$. Therefore the expressions $E=x_{\sigma_1}L_1x_{\sigma_2}L_2\ldots L_{k-1}x_{\sigma_p}$ are in bijective correspondance to  maps $f_E$ given by
$f(\sigma_1)=\ldots =f(\sigma_{a_1})=1$, $f(\sigma_{a_1+1})=\ldots =f(\sigma_{a_2})=2$,\ldots, 
$f(\sigma_{a_{p-\ell-2}+1})=\ldots =f(\sigma_{a_{p-\ell-1}})=p-\ell-1$, and $f(\sigma_{a_{p-\ell-1}+1})=\ldots =f(\sigma_p)=p-\ell$.

According to (6.8) of \cite{dw}, the number of surjections of a set of $m$ elements onto a set of $n$ elements equals 
$n!S(m,n)$, where 
$S(m,n)$ is the Stirling number of the second kind. 
Thus $c_{p,\ell}=(p-\ell)!S(p,p-\ell)$.

If we substitute $x=-n$ in the formula
$x^p=\sum_{j=1}^r S(p,j) x(x-1)\ldots (x-j+1)$
from Theorem 6.10 of \cite{dw}, we obtain 
\[n^p=\sum_{j=1}^p (-1)^{p-j} j!S(p,j) F^j_n=\sum_{\ell=0}^{p-1} (-1)^{\ell} c_{p,\ell} F^{p-\ell}_n.\quad \qed\]

Note that an earlier (and much longer) proof of this statement was given in \cite{cer}.


\end{document}